\documentclass[12pt]{amsart}
\usepackage{amssymb}
\usepackage{graphicx}

\theoremstyle{plain}
\newtheorem{theorem}{Theorem}[section]
\newtheorem{lemma}[theorem]{Lemma}

\newtheorem{corollary}[theorem]{Corollary}
\theoremstyle{definition}
\newtheorem{definition}[theorem]{Definition}
\newtheorem{remark}[theorem]{Remark}

\numberwithin{equation}{section}

\renewcommand{\phi}{\varphi}

\begin{document}

\title{The Combinatorics of Avalanche Dynamics}
\author{Manfred Denker}
\address{Mathematics Department, Pennsylvania State University, State College PA 16802, USA} \email{denker@math.psu.edu}
\author{Ana Rodrigues}
\address{Matematiska institutionen, KTH, Lindstedtsv\"agen 25, S-100 44
Stockholm, Sweden.} \email{ana.rodrigues@math.kth.se}
\thanks{The research of MD was supported by the National Science Foundation grant DMS-1008538. The research of AR is supported by the Swedish Research Council (VR Grant 2010--5905). The second author would like to thank the G\"oran
Gustafsson Foundation UU/KTH for financial support.}
\date{\today}

\maketitle

 \begin{abstract} We give a simple and elementary proof of  the identity
$$\sum_{r=1}^n\sum_{k_1,...,k_r\ge 1: \sum_{i=1}^r k_i= n} \frac {n!}
{k_1!k_2!...k_r!}k_1^{k_2}...k_{r-1}^{k_r}=(n+1)^{n-1}$$
where $n\in \mathbb N$. A first application of this formula shows
Cayley's theorem \cite{Caley} on the number of
trees with $n+1$ vertices (in fact the formula is equivalent to Cayley's result).  A second application gives the distribution of avalanche sizes, which can be deduced for
general dynamical systems and also as a bilogically motivated urn model in
probability. In particular,  the law 
 of avalanche sizes in Eurich et al. \cite{EHE} and
Levina \cite{Levina} is closely related to this dynamical representation.
\end{abstract}

\section{Introduction}

The purpose of this note is to introduce avalanche dynamics into the theory of dynamical systems. This is a type of inducing scheme which defines new dynamical systems leaving the orbit structure mostly untouched. Such inducing schemes are first of all known as inducing on a set: let $(X,T)$ be a dynamical system (discrete time) and $A\subset X$, then the transformation induced by $A$ is defined by
$$ T_A(x) = T^{\phi(x)}(x)$$
where $\phi(x) = \inf\{n\ge 1: T^n(x)\in A\}$. This defines a transformation on all points in $A$ which return to $A$ infinitely often. Next we mention Schweiger's jump transformation (see \cite{Aa}) which is defined by
$$ T^*_A(x) = T^{\psi(x)}(x)$$
where $\psi(x)=1$ for $x\not \in A$ and  $\psi(x) = \inf \{n\ge 1: T^n(x)\not\in A\}$ for $x\in A$. Variants of this are also known, like Young tower constructions (see \cite{Young}).

Avalanche dynamics is another inducing scheme and will be defined in Section 3. The new orbits will be a subset of the old orbits leaving the order structure invariant. This definition is motivated by the study of avalanches in neural dynamics (see eg \cite{Levina}).

We begin with a short and rough description of the physiological background of neuronal dynamics. There are two types of cells in the central nervous system, one of them is called neuron, which communicate by sending and receiving electrical impulses. The cell membrane has built into it channels and ion pumps, letting potassium ions rushing out, and sodium ions flushing in. This process of exchange of potassium and sodium ions (once initiated) stops after a few milliseconds when repolarisation is achieved. The process needs an external activation for getting started.

 The integrate and fire model goes back to Lapique in 1907 (\cite{Lapique}). It describes the time series of the potential of the cell membrane. It can be written as an ordinary differential equation of the form
$$ C_m \frac{dV(t)}{dt}= g_l \left(V_{res}-V(t)\right) +I(t),$$
where $C_m$ is called the capacitance of the cell membrane, $V(t)$
is the voltage (potential) at time $t$, $V_{res}$ is the residual potential, $g_l$ is the leak conductivity of the cell membrane and $I(t)$ is the current at time $t$. This time evolution of the potential of a neuron is interrupted if the potential reaches a certain threshold value (this is when the process of exchange of potassium and sodium ions starts). At this point it gives a fixed electrical impulse to each of its neighboring neurons.  Thus other neurons' potential may reach the threshold, thus initiating an avalanche of firing neurons. We are concerned in Section 3 with the size of this avalanche. If the dynamics of the neurons is stationary, we derive an asymptotic expression for the probability of the avalanche size. Such expression had been found earlier by Eurich et al (\cite{EHE}) and Levina (\cite{Levina}). Our expression is slightly more general and includes the earlier results as conditional distributions given that one neuron is excited. Furthermore, we also study this distribution in more detail.

The key observation to our approach is a combinatorial formula which we prove
in Section 2. Although this formula does not seem to be discussed in the
literature, it can be deduced from Cayley's result on the number of labeled
trees (\cite{Caley}); it also is a special case of Corollary 8 in
\cite{Pit}. In fact we show that the formula also implies Cayley's result,
thus adding a new proof of Cayley's formula.  In view the combinatorial formula and the two applications, we suspect that there is a connection between avalanche sizes and branching processes in the critical case (avalanche sizes are calculated by one side in the formula, while the total number of successors in a branching process are calculated from a rooted tree, which is represented by the other side of the formula). In fact, Levina has derived such a relation showing that the asymptotic distributions are the same. This leads to a power law of the distribution using a  result by Otter (\cite{Otter}). Such a power law is also derived in Section 3.

{\it Acknowledgement:} We would like to thank Ira Gessel and Wlodek Bryc for some helpful remarks.

\section{A Combinatorial Formula} In 1889 Cayley showed that the number of
labeled trees  with $n$ distinguishable vertices is $n^{n-2}$. We begin giving an apparently new proof of this result. The article of Moon (\cite{Moon}) lists ten different proofs. Renyi (\cite{Ren}) gave another proof of this  fact. As noticed in the beginning of the proof by Clarke (\cite{Clarke}), all trees labeled with $n+1$ points can be represented as rooted trees where an arbitrary chosen vertex is fixed as the root. Counting these
trees by dividing the remaining $n$ vertices into $r$ subsets $V_1$,...,$V_r$ of sizes $k_1$,...,$k_r$ and letting $V_l$ denote the vertices at distance $l$ from the root (in the path lengths metric), one has $k_{l-1}^{k_l}$ choices to connect to the set $V_{l-1}$. Given the Cayley result we thus have a proof of the following Theorem  \ref{combi}. Conversely, we shall prove Theorem \ref{combi} by simple induction, which leads to a new proof of  Cayley's Theorem as an obvious corollary.

\begin{theorem}\label{combi} Let $n\in \mathbb N$, the natural numbers. Then
$$\sum_{r=1}^n\sum_{k_1,...,k_r\in \mathbb N: \sum_{i=1}^r k_i= n} \frac {n!}
{k_1!k_2!...k_r!}k_1^{k_2}...k_{r-1}^{k_r}=(n+1)^{n-1}.$$
\end{theorem}

\noindent{\it Proof.}  The binomial formula reads as
\begin{equation*}\label{basic}
n^{n-k-1}= (n-k+k)^{n-k-1}
= \sum_{j=1}^{n-k} {n-k\choose j} jk^{j-1}(n-k)^{n-k-j-1}.
\end{equation*}
Now proceed by induction to show for $s\ge 1$ that
\begin{eqnarray*}
(n+1)^{n-1}&=&\sum_{r=1}^s\sum_{k_1,...,k_r\ge 1;k_1+...+k_r=n} \frac{n!}{k_1!...k_r!} k_1^{k_2}k_2^{k_3}...k_{r-1}^{k_r}\\
&&+
\sum_{k_1,...,k_s\ge 1: k_1+...+k_s<n} \frac{n!}{k_1!...k_s!(n-k_1-...-k_s)!}\\
&&\quad k_sk_1^{k_2}...k_{s-1}^{k_s}(n-k_1-...-k_{s-1})^{n-k_1-k_2-...-k_s-1}.
\end{eqnarray*}

Having established the basic formula we immediately derive the following
corollary.

\begin{corollary} Let $n\ge 1$. Then
$$\sum_{r=1}^n\sum_{k_1,...,k_r\ge 1: \sum_{i=1}^r k_i= n} \frac {n!}
{k_1!k_2!...k_r!}k_1^{k_1-1}...k_{r}^{k_r-1}=(n+1)^{n-1}.$$
\end{corollary}

\noindent{\it Proof.} As mentioned before the previous theorem every tree with $n+1$
labeled vertices  $\{v_0,...,v_n\}$ can be represented as a rooted tree with (say) root $v_0$.  Consider a partition of $\{v_1,...,v_n\}$ into subsets $E_1$,...,$E_r$ of cardinalities $k_1$,...,$k_r$. The number of labeled trees with vertices from $E_l$ is $k_l^{k_l-2}$ according to Cayley's formula. There are $k_l$ choices for a root in $E_l$. For each choice of a root in $E_l$ and each tree in $E_l$, $l=1,...,r$,  we can construct a unique tree of all vertices by connecting the roots in $E_l$ with $v_0$.
This has $\prod_{l=1}^r k_l^{k_l-1}$ choices. Summing over $E_1$,...,$E_r$, then over $k_1,...,k_r\ge 1$ with $k_1+...k_r=n$, and and finally over $r=1,...,n$ shows the corollary.

\section{Abelian Distributions}

In 2002 Eurich et al. \cite{EHE} proposed a probability which
describes the sizes of avalanches in neural dynamics. The proportionality
factor in this approach has been determined  by
Levina \cite{Levina} in 2008. In addition she was able to determine the
expectation of this distribution.

 Suppose we have $N$ neurons, of which one is in its
firing state. Then it gives a certain internal impulse to each other neuron
(a complete network is assumed). There may be a certain number of inactive
neurons which move to the firing state after receiving the internal
input. They begin to fire sending the same impulse to all other
neurons. Continuing this way one obtains an avalanche, and if it stops, denote
by $X$ the total number of firing neurons in an avalanche period. It was claimed in \cite{EHE} that
the distribution of $X$ has asymptotically (as the number of neurons $N$ tends to $\infty$ and the internal impuls $Np \to 1$) the following form:

\begin{theorem}\label{levina} (\cite{Levina}) Let $N$ be a positive integer and $p\in [0,\frac {1}{N})$. Then for $k=1,2,...,N$
\begin{equation}\label{abelian} p_k = P(X=k)=\frac{1-Np}{1-(N-1)p}  k^{k-2} {N-1\choose k-1}p^{k-1}(1-kp)^{N-k-1}
\end{equation}
defines a probability distribition with expectation
\begin{equation}\label{expectation}  E[X]= \frac 1{ 1-(N-1)p}.
\end{equation}
\end{theorem}

This distribution (\ref{abelian}) was called Abelian distribution in \cite{Levina}.
The motivation of this formula is as follows. Suppose one neuron is excited in
a complete neural network. Then $k-1$ is the number of neurons firing
successively in the avalanche, each with probability $p$. Note that $(1-(kp)$
is the probability that a neuron is not firing during the avalanche period of avalanche size $k$. In fact, we would expect
(and derive it below) that the exponent of the factor $1-kp$ is
$N-k$ the number of non-firing neurons.

There are two proofs of Theorem \ref{levina}, one is given by Levina by showing that
numerator and denominator are polynomials in $p$ and then verifying the identity
at a finite number of points. The second method is by using a generalized
binomial theorem (personal communication by M. Rao to one of us). Below we will give a dynamical
proof of a related formula using measure preserving dynamical systems. Besides we will show how
such a  law can be derived  using Theorem \ref{combi}.
\\

\centerline{\bf Urn models}

Consider $N$ distinguishable balls and $M$ enumerated urns $1,2,...,M$. Distribute the $N$ balls randomly in the urns, so the probability of an elementary event is $\frac 1{M^N}$.

Let $X$ be the largest integer $r\in \{1,...,M\}$ such that for each $k\le r$ the number of balls in urns $1,...,k$ is at least $k$. For example, $X=0$ if no balls is placed in urn $1$, and $X=1$ if one ball is placed into urn $1$ and no ball in urn $2$. In general, one has $r$ balls in the urns numbered $1,...,r$  and none in the $r+1$st urn and the number of balls in urns $1,..,k$ is at least $k$ for each $k\le r$. We have

\begin{lemma} The distribution of $X$ is given by the avalanche distribution:
$$P(\{X=a\})= {N\choose a}(a+1)^{a-1}\frac 1{M^a}\left(1-\frac{a+1}M\right)^{N-a},$$
where $a=0,1,2,...,N$.
\end{lemma}

\noindent{\it Proof.} Partitioning $\{1,2,...,N\}$ into sets $J_1,...,J_{r+1}$ of cardinalities $k_1$, $k_2$,..., $k_{r+1}=N-a$
the number of elementary configurations for which $X=a=k_1+...+k_r$ is
$$ k_1^{k_2}k_2^{k_3}...k_{r-1}^{k_r}.$$
(Put the balls from $J_1$ into urn 1, then the balls in $J_2$ in one of the boxes $2,...,k_1+1$  etc.)
If $a$ is fixed, sum over all such partitions to arrive at
$$ P(X=a)= \sum_{r=0}^N \sum_{k_1,...,k_r\ge 1; k_1+...k_r=a} \frac{a!}{k_1!...k_r!} k_1^{k_2}...k_{r-1}^{k_r}  {N\choose a} \frac 1{M^a} (1-\frac{a+1}M)^{N-a}.$$
Now apply Theorem \ref{combi} to obtain the formula in the lemma.\\

\centerline{\bf Dynamical systems as a model}

\noindent For $i=1,...,N$ let $S_i:X_i\to X_i$ be continuous maps on the metric spaces $X_i$ with metric
$d_i$ for $i=1,...,N$.
Let $X=X_1\times X_2\times ...\times X_N$ denote the $N$-fold product space and
$S=S_1\times S_2\times...\times S_N$ the product map on $X$.

We define the avalanche associated with open sets $U_i\subset
X_i$ ($1\le i\le N$) which satisfy the condition that that they are the top level of towers
of heights $n_i+1>N$ for the transformation $S_i$. That is:
there exist sets $B_i\subset X_i$ such that the sets $S_i^k(B_i)$ for
$k=0,...,n_i$ are pairwise disjoint and $S_i^{n_i}(B_i)=U_i$.
An avalanche of $(X,S)$ at the point $x\in X$ is defined by the set of
coordinates which move under $S_i^j$ ($j\ge 0$) through the sets $U_i$ such that at each
level $n_i-l$ of the tower there are at least $l$ coordinates which are above
that level. A coordinate in a set $U_i$ or moving through $U_i$ during an
avalanche period is called an excited coordinate (state).

 Formally the avalanche size is defined as follows. For $x=(x_1,x_2,...,x_N)\in X$ define
$$ A(x,1)=|\{1\le i\le N: x_i\in U_i\}|$$
and
$$ A(x,2)= |\{1\le i\le N: \exists\ 0\le l\le A(x,1)\ni S_i^l(x_i)\in U_i\}|.$$
Then, recursively, define for $k\ge 2$
$$ A(x,k+1)= |\{1\le i\le N: \exists 0\ \le l\le A(x,k)\ni S_i^l(x_i)\in
U_i\}|.$$
Note that the sequence $A(x,k)$ is increasing.\\

\begin{lemma} There exists a minimal $k\le N$ such that $A(x,k+1)=A(x,k)$.
\end{lemma}

\noindent{\it Proof.} Since the iterates $S_i^l(B_i)$ of $B_i$ are pairwise
disjoint for $l=0,...,N<n_i$  the i-th coordinate can fall only at most once
into  $U_i$
for iterations up to time  $N$. \\

\begin{definition} The value $A(x)=\mbox{card} A(x,k)$ such that
$$ A(x,k)=A(x,k+1)$$
is called the avalanche size at $x\in X$.
\end{definition}

\noindent Suppose we have a partition of $\{1,2,...,N\}$ into $r+1$ sets of sizes
$k_1,...,k_{r+1}$, so $k_1+k_2+...+k_{r+1}=N$. We denote the sets of the partition by
$I_1,...,I_{r+1}$. Such a partition corresponds  to an
avalanche of size $a=k_1+k_2+...+k_r$: letting the coordinates (states) in $I_1$ be excited in the
beginning, then the coordinates in $I_2$ until we let  the coordinates
in $I_r$ be excited; finally we denote
by $I_{r+1}$ the coordinates which are never excited.

\noindent We describe this formally. Let
$$ E_1=\bigotimes_{i\in I_1} S_i^{n_i}(B_i)\otimes\bigotimes_{i\not\in I_1}X_i,$$
so for $x\in E_1$ the coordinates $x_i$ of $x$ for which $i\in I_1$ belong to the
sets $U_i=S^{n_i}(B_i)$, and are excited states.

\noindent In the second step we have the particles in $I_2$ firing if $x$
belongs to
$$ E_2=\bigotimes_{i\in I_2} \bigcup_{l=1}^{k_1}
S_i^{n_i-l}(B_i)\otimes\bigotimes_{i\not\in I_2} X_i.$$

\noindent Continuing this way we obtain sets
$$E_l= \bigotimes_{i\in I_l} \bigcup_{j=1}^{k_{l-1}}
S_i^{n_i-k_1-...-k_{l-2}-j}(B_i)\otimes\bigotimes_{i\not\in I_l} X_i,$$
which describes the points for which the coordinates are firing in the l-th step but not before ($l=2,...,r$).

\noindent Finally we put
$$E_{r+1}= \bigotimes_{i\in I_{r+1}} X_i\setminus \bigcup_{l=0}^{a}
S_i^{n_i-l}(B_i)\otimes\bigotimes_{i\not\in I_{r+1}} X_i, $$
where $a=k_1+...+k_r$ is the total number of firing coordinates. The set $E_{r+1}$ has only coordinates which are not at all
firing. Also note that $A(x)=a$ is the avalanche size.

\noindent Now put $E=E_1\cap E_2\cap...\cap E_{r+1}$ and this set is where we have
an avalanche of size $a$ with exited states described by the sets $I_i$, $i=1,...,r$.\\

\noindent We assume now that there exist  $S_i$-invariant ergodic measures $m_i$ on $X_i$, $1\le i\le N$. Then the
product measure
$$ m=\otimes_{i=1}^N m_i$$
is $S$-invariant and ergodic. We have

\begin{lemma} Let $I_1$, ..., $I_{r+1}$ be as above. Then
$$m(E)= \prod_{i\in I_1} m_i(U_i) \prod_{l=2}^r \prod_{i\in I_l}
  k_{l-1}m_i(U_i)\prod_{i\in I_{r+1}} (1-(a+1)m_i(U_i)).$$
\end{lemma}

\begin{corollary} let $m_i(U_i)=p$ for all $i=1,...,N$. Then
$$ m(E)=p^a(1-(a+1)p)^{N-a}\prod_{l=1}^r k_{l-1}^{k_l}.$$
\end{corollary}

\begin{theorem}\label{avalanche} Let $A:X\to \mathbb N$ denote the
avalanche size function. Define
$$K_r(a)=\{(k_1,k_2,...,k_{r+1}): k_1,...,k_r\ge 1, k_{r+1}=N-a;\ \sum_{i=1}^r
k_i=a\}.$$
Then
for any $a=0,1,2,...,N$
\begin{eqnarray*}
&& m(\{x\in X: A(x)=a\})\\
&=&\sum_{r=1}^a \sum_{(k_1,k_2,...,k_{r+1})\in K_r(a)}
\sum_{I_1,...,I_{r+1}}  \prod_{i\in I_1} m_i(U_i) \prod_{l=2}^r \prod_{i\in I_l}
  k_{l-1}m_i(U_i)\\
&& \prod_{i\in I_{r-1}} (1-(a+1)m_i(U_i)),
\end{eqnarray*}
where $\sum_{I_1,...,I_{r+1}}$ denotes summation over all partitions of
$\{1,...,N\}$ into sets $I_1,...,I_{r+1}$ of sizes $k_1,...,k_{r+1}$.

In particular, if $m_i(U_i)=p$ for all $i=1,...,N$ then
$$m(\{x\in X: A(x)=a\})= (a+1)^{a-1} {N\choose a}p^a (1-(a+1)p)^{N-a}.$$
\end{theorem}

\noindent{\it Proof.} Note that the first formula is immediate from the
foregoing discussion. Moreover, if all $m_i(U_i)=p$ for some $p\in [0,1]$ then, by Theorem \ref{combi},
the formula reduces to
\begin{eqnarray*} && m(\{x\in X: A(x)=a\})\\
&=&\sum_{r=1}^a \sum_{k_1+k_2+...+k_r=a;k_{r+1}=N-a}
\frac{N!}{k_1!k_2!...k_{r+1}!}  p^a \prod_{l=2}^{r} (k_{l-1})^{k_l}
(1-(a+1)p)^{N-a}\\
&=& \sum_{r=1}^N \sum_{k_1+k_2+...+k_r=a}
\frac{a!}{k_1!k_2!...k_{r}!} \prod_{l=2}^{r} (k_{l-1})^{k_l}
{N\choose a} p^a  (1-(a+1)p)^{N-a}\\
&=& (a+1)^{a-1}{N\choose a} p^a (1-(a+1)p)^{N-a}.
\end{eqnarray*}

\noindent{\bf Corollary:} Let $p=\alpha/N$. Then for each $a$
$$ \lim_{\alpha\to 1}\lim_{N\to\infty} \log \frac{m(\{A=a\})}{m(\{A=a+1\}}= -\frac 3{2a} +o(\frac 1 a).$$

\noindent{\it Proof.} The proof uses Taylor expansion of the logarithm:
\begin{eqnarray*}
\lim_{\alpha\to 1}\lim_{N\to\infty}\log \frac{m(\{A=a\})}{m(\{A=a+1\}} &= & 1+ a\log \frac {a+1}{a+2}\\
&=&  1+ a(-\frac{1}{a+2}- \frac{1}{2(a+2)^2}) + o(1/a) \\
&=&\frac{2}{a+2}- \frac{a}{2(a+2)^2}+o(1/a)\\
&=& \frac 3{2a} +o(1/a).
\end{eqnarray*}

\begin{remark} Note that this asymptotic means that
$$ m(x\in X: A(x)=a)\sim a^{-3/2}.$$
\end{remark}

\begin{remark} The distribution has only local maxima for $\alpha$ close to 1, i.e. in the supercritical case.
\end{remark}

\begin{remark} It is straight forward to deduce a  Levina type result from Theorem
  \ref{avalanche}. Since in neural networks it is assumed that one neuron
  starts to fire, we are looking at the conditional distribution that one
  particular neuron is firing. Conditioned on this event, there are $N-1$
  neurons remaining which may form an avalange of possible sizes $1,...,N$
  including the initial firing neuron. According to Theorem \ref{avalanche} the
  distribution is given by
$$P(A=a)= a^{a-2}{N-1\choose a-1} p^{a-1} (1-ap)^{N-a},$$
where $a=1,...,N$.
This is almost Levina's formula; the power in the last factor differs.
\end{remark}

We can get the expectation of the Abelian distribution in Theorem \ref{levina}
from our theorem.

\begin{corollary} We have
$$  \frac{1-Np}{1-(N-1)p}\sum_{k=1}^N k^{k-1}{N-1\choose k-1} p^{k-1}(1-kp)^{N-k-1}= \frac {1}{1-(N-1)p}.$$
\end{corollary}

\begin{remark} The corollary means that the avalanche size
  in Theorem \ref{levina} has expectation $\frac {1}{1-(N-1)p}$ (see (\ref{expectation})).
\end{remark}

\noindent{\it Proof.} For each $x\in [0,\frac 1N)$ we know from Theorem \ref{avalanche} that
$$ f(x)=\sum_{a=0}^N (a+1)^{a-1} {N\choose a}x^a (1-(a+1)x)^{N-a}=1.$$
Taking derivative with respect to $x$ yields
\begin{eqnarray*}
 f'(x)&=& (N+1)^{N-1}Nx^{N-1}- N(1-x)^{N-1} \\
&+&\sum_{a=1}^{N-1}
(a+1)^{a-1}{N\choose a} x^{a-1}(1-(a+1)x)^{N-a-1}(a(1-Nx)-Nx) =0.
\end{eqnarray*}
Since $N= \frac{N(1-Nx)-Nx}{1-(N+1)x}$ we obtain
$$ f'(x)= \sum_{a=0}^{N}
(a+1)^{a-1}{N\choose a} x^{a-1}(1-(a+1)x)^{N-a-1}(a(1-Nx)-Nx) =0.$$
Multiplying by $x$ and replacing $a$ by $b-1$ yields
\begin{equation}\label{identity} f'(x)= \sum_{b=1}^{N+1}
b^{b-2}{N+1-1\choose b-1} x^{b-1}(1-bx)^{N+1-b-1}((b-1)(1-Nx)-Nx) =0.
\end{equation}
By Theorem \ref{levina} we get
$$\frac{1-(N+1)x}{1-Nx}\sum_{b=1}^{N+1}
b^{b-2}{N\choose b-1} x^{b-1}(1-bx)^{N-b}=1.$$
Thus, it follows from (\ref{identity}) that
\begin{eqnarray*} && \frac{1-Nx}{1-(N-1)x}\sum_{b=1}^{N+1}
b^{b-1}{N\choose b-1} x^{b-1}(1-bx)^{N-b}(1-Nx)\\
&&=\frac{1-Nx}{1-(N-1)x}\sum_{b=1}^{N+1}
b^{b-2}{N\choose b-1} x^{b-1}(1-bx)^{N-b}\\
&&=1.
\end{eqnarray*}

\end{document}